\documentclass[12pt]{article}
\usepackage{amssymb}
\usepackage{latexsym}

\newtheorem{Theorem}{Theorem}[part]

\newtheorem{Proposition}{Proposition}[part]

\newtheorem{Lemma}{Lemma}[part]

\newtheorem{Remark}{Remark}[part]

\makeatletter \@addtoreset{equation}{section}

\@addtoreset{Definition}{section}

\@addtoreset{Theorem}{section}

\@addtoreset{Proposition}{section}

\@addtoreset{Property}{section}

\@addtoreset{Assumption}{section}

\@addtoreset{Corollary}{section}

\@addtoreset{Lemma}{section}

\@addtoreset{Remark}{section}

\@addtoreset{Example}{section}

\def \proof{{\noindent \bf Proof. }}
\def \ep{\hbox{ }\hfill$\Box$}
\def\reff#1{{\rm(\ref{#1})}}

\def\Ac{{\cal A}}

\def\Dc{{\cal D}}

\def\Fc{{\cal F}}
\def\Gc{{\cal G}}
\def\Hc{{\cal H}}

\def\Zc{{\cal Z}}

\def\Zc{{\cal Z}}

\def\no{\noindent}
\def\Pas{\mathbb{P}-\mbox{a.s.}}
\def\Fbb{{\rm \bf F}}
\def\x{\times}

\def\05{\frac{1}{2}}
\def\-1{^{-1}}
\def\1{{1\hspace{-1mm}{\rm I}}}
\def\={\;=\;}
\def\.{\;.}

\def\be{\begin{eqnarray}}
\def\ee{\end{eqnarray}}

\def\b*{\begin{eqnarray*}}
\def\e*{\end{eqnarray*}}

\def\And{\;\mbox{ and }\;}

\def\Esp#1{\mathbb{E}\left[#1\right]}

\def\Pro#1{\mathbb{P}\left[{#1}\right]}

\def\And{\mbox{ and } }
\def\pourtout{\mbox{ for all } }

\def \E{\mathbb{E}}
\def \F{\mathbb{F}}
\def \H{\mathbb{H}}
\def \M{\mathbb{M}}
\def \N{\mathbb{N}}
\def \R{\mathbb{R}}

\def\P{\mathbb{P}}
\def\Q{\mathbb{Q}}

\def\T{\mathbb{T}}
\def\L{\mathbb{L}}

\def\ri#1{\mbox{\rm ri}(#1)}

\def\vs#1{\vspace{#1mm}}

\def\i{{\rm (i) }}
\def\ii{{\rm (ii) }}
\def\iii{{\rm (iii) }}

\def\HKP{\mbox{\rm{\bf KP}}}
\def\HNF{\mbox{\rm{\bf HN$^0$}}}
\def\HFa{\mbox{\rm{\bf HF$_1$}}}
\def\HFb{\mbox{\rm{\bf HF$_2$}}}
\def\NAw{\mbox{\rm{\bf NA$^w$}}}
\def\NAs{\mbox{\rm{\bf NA$^s$}}}
\def\NAr{\mbox{\rm{\bf NA$^r$}}}
\def\HL{\mbox{\rm{\bf H$\Ac$}}}
\def\DFa{\mbox{\rm {\bf D$_1$}}}
\def\DFb{\mbox{\rm {\bf D$_2$}}}
\def\EF{\mbox{\rm {\bf EF}}}

\begin{document}

\title{No-arbitrage in discrete-time   markets with proportional transaction costs and general information structure}
\date{December 28, 2004}
\vs5

\author{Bruno BOUCHARD\footnote{Web page:
http://felix.proba.jussieu.fr/pageperso/bouchard/boucharda.htm} \footnote{I am grateful to Fabian Astic for  his remarks.} 
                        \\ \small Laboratoire de Probabilit{\'e}s et
              Mod{\`e}les Al{\'e}atoires
             \\ \small CNRS, UMR 7599,
              Universit{\'e} Paris 6
             \\ \small  and CREST
             \\ \small email: bouchard@ccr.jussieu.fr
             }
\maketitle




\begin{abstract}
We discuss the no-arbitrage conditions in a general framework
 for  discrete-time models of  financial markets with
proportional transaction costs and general information structure.
We extend the results of Kabanov and al. (2002),
 Kabanov and al. (2003) and Schachermayer (2004) to the case where bid-ask spreads
are not known with certainty. In the ``no-friction" case, we
retrieve the result of Kabanov and Stricker (2003).
\end{abstract}

\vspace{5mm}

\no Key words: Absence of arbitrage, proportional transaction
costs, imperfect information, optional projection.

\vspace{2mm}

\no JEL Classification Numbers: G10.

\vs1

\no AMS (2000) Subject Classification: 91B28, 60G42.

\section{Introduction}

While the ``insider trading" problem, where the agent's filtration
$\H$ is strictly bigger than the asset's filtration $\F^S$, has
been widely studied in the recent literature, see e.g.
\cite{ABS}, \cite{CIKN}, \cite{I} and the references therein, less
care has been given to the imperfect information case where $\H$
does not contain $\F^S$.  Such situations may arise for instance
if the small investor has not a direct access to the market. In
this case, his orders can be executed with a delay and therefore
at a price which is not known in advance, see \cite{Gaussel} and
\cite{KSdelayed}. From the point of view of the arbitrage
theory,   the ``insider trading" case is well known. Indeed, all
the necessary and sufficient condition for the absence of
arbitrage opportunities available in the case $\H=\F^S$ apply to
the general case $\F^S\subset \H$. In particular, the usual
``no-arbitrage" conditions imply that, roughly speaking (see
\cite{DSborne} and \cite{DSnonborne} for precise results), prices
must be semi-martingales in the filtration $\H$ and   that there
must be an equivalent probability measure $\Q$ under which they
are $(\Q,\H)$-local martingales.

The case of imperfect information   where $\H \subset \F^S$ and
$\H\neq \F^S$    is much more difficult to handle.  In particular,
the arguments of \cite{DSborne} do not work in this situation.
Even in the case of infinite discrete time, the proof of
\cite{Sinfinite} does not apply. However, in finite discrete time,
it was noticed in \cite{KSdelayed} that the proof of the
Dallang-Morton-Willinger theorem reported in \cite{teachers note}
still holds up to minor modifications for any given filtration
$\H$. In this case, the no-arbitrage condition is equivalent to
the existence of a probability measure $\Q$ such that the optional
projection under $\Q$, $(\E^{\Q}[S_t~|~\Hc_t])_t$, of the asset
prices $(S_t)_t$ on $\H=(\Hc_t)_t$ is a $(\Q,\H)$-martingale.

The aim of  this paper is to extend this result to the case where
exchanges are subject to proportional transaction costs. In the
recent literature, such models have been widely studied, from the
seminal work of \cite{JKequi} to the recent papers
\cite{KabStri}, \cite{KSR01}, \cite{KSR04}, \cite{schach},
\cite{rasThese}, \cite{Penner} and \cite{BP04} among others. The recent
abstract formulation consists in introducing a sequence of random
closed convex cones $(K_t)_t$ and describing the wealth process as
$V_t=\sum_{s\le t} \xi_s$ with $\xi_s \in -K_s$ a.s. The ``usual"
example is given by
    \be\label{eq modele usuel}
    -K_t(\omega)= \{x\in \R^d~:~\exists\;a\in \M^{d}_+,\;
    x^i\le \sum_{j\le d} a^{ji}-a^{ij}\pi^{ij}_t(\omega)\;,\;i\le d \}\;,
    \ee
where $\M^d_+$ denotes the set of square $d$-dimensional matrices
with non-negative entries. Here $\pi^{ij}$ should be interpreted
as the costs in units of asset $i$ one has to pay to obtain one
unit of asset $j$. If we allow to throw out money, an exchange
$\xi_t$ at time $t$ is then affordable if $\xi_t \in -K_t$ a.s.

In the case of imperfect information, i.e. $\pi$ is not
$\H$-adapted, this approach cannot be used since $K$ is no longer
$\H$-adapted. Hence, we have to change the modelisation. Instead
of considering the $\xi$'s as the controls, we have to rewrite
them as $\xi^i_t$ $=$ $\sum_{j\le d}
\eta_t^{ji}-\eta_t^{ij}\pi^{ij}_t$, where $\eta$ is an
$\H$-adapted process with values in the set of square
$d$-dimensional matrices with non-negative entries $\M^{d}_+$.
Here, $\eta_t^{ji}$ is the number of physical units of $i$ we
 obtain, at time $t$, against $\eta_t^{ji}\pi^{ji}_t$ units  of
$j$. Because the $\xi$ may not be adapted the proofs of
\cite{KSR01}, \cite{KSR04} and \cite{schach} does not apply to
this setting and, in contrast to \cite{KSdelayed}, we have to work a
bit more to extend their results.

In the above model, we fix the number of units $\eta_t^{ij}$ of
asset $j$ we want to buy and the number of units of asset $i$ one
has to sell is given by  $\eta_t^{ij}\pi^{ij}_t$. In the case of
perfect information, i.e. $\pi$ is $\H$-adapted, one can also
fix  the amount $\tilde \eta_t^{ij}$ of units of asset $i$
one wants to sell and compute $\eta_t^{ij}$ accordingly by using
the formula $\tilde \eta_t^{ij}=\eta_t^{ij}\pi^{ij}_t$. But in
the case where $\pi$ is not $\H$-adapted this is no more possible
and one cannot control exactly $\tilde \eta_t^{ij}$. This means
that orders can be formulated only in terms of the quantity  of
units of the asset we want to buy and we shall see  in Subsection
\ref{subsec currency 2} that, in such a situation, orders may be
non-reversible even in the case of no-friction where
$\pi^{ij}=\pi^{ji}$ for all $i,j\le d$. Clearly, this is not
reasonable and in practice one should also be able to fix $\tilde
\eta_t^{ij}$. To pertain for such orders, one can slightly modify
the above model by taking $\xi$ in the form $\xi^i_t$ $=$
$\sum_{j\le d}
\eta_t^{ji}(1+\lambda^{ij}_t(\omega)\1_{\eta_t^{ji}<0})
    -\eta_t^{ij}\tau^{ij}_t(\omega)(1+\lambda^{ij}_t(\omega)\1_{\eta_t^{ij}>0})$ where $\eta$ is $\H$-adapted process
with values in  the set $\M^d$ of square $d$-dimensional
matrices. Here, $\tau^{ij}$ stands for the costs in units of asset $i$ one
has to pay to obtain one unit of asset $j$, before to pay the
transaction costs. The transaction costs
$\eta_t^{ij}\tau^{ij}_t\lambda^{ij}(\omega)\1_{\eta_t^{ij}>0}$ $+$
$\eta_t^{ji}  \lambda^{ji}_t(\omega)\1_{\eta_t^{ji}<0}$ are paid in units of the sold
asset $i$. With the above notation, one has
$\pi^{ij}_t=\tau^{ij}_t(1+\lambda^{ij})$. This corresponds to
 \be\label{eq modele 1}
    -K_t(\omega)&=& \{x\in \R^d~:~\exists\;a\in \M^{d},\;
    \\
    &&
    x^i
    \le
    \sum_{j\le d}
    a^{ji}(1+\lambda^{ij}_t(\omega)\1_{a^{ij}<0})
    -a^{ij}\tau^{ij}_t(\omega)(1+\lambda^{ij}_t(\omega)\1_{a^{ij}>0})
    \;,\;i\le d \}\;,\nonumber
 \ee
Contrary to the model
\reff{eq modele usuel}, we can now fix the number of units of
asset $i$ we want to sell against units of asset $j$ by fixing
$\eta_t^{ji}<0$ so that $|\eta_t^{ji}|$ coincides with the amount
of exchanged units of $i$. Once again, in the case of perfect
information both models are equivalent, but this is no more
true   if $\tau$ and/or $\lambda$ are not $\H$-adapted. One could
also argue that paying the transaction costs in units of the
asset which is sold, as in \reff{eq modele 1}, is not the same
thing than paying these costs in units of the asset which is
bought. Here  again one could consider a more general model which
pertains for different costs structures.

In order to take into account all these different situations, we
propose a general formalism where the wealth process $V$ is
written as $V_t=\sum_{s\le t} F_s(\eta_s)$, for some sequence of
random maps $F$ $=$ $(F_t)_t$.  Here,  $\eta$ is $\H$-adapted
process with values in a closed convex cone $\Ac$ of  $\M^{d}$ (we
have in mind to take $\Ac$ $=$ $\M^{d}$, however,  in order to
take also the model \reff{eq modele usuel} into account it is
more convenient to allow for the possibility of having $\Ac$ $=$
$\M^{d}_+$). We make no assumption on the filtration under which
$F$ is adapted. Thus, this approach pertains for the cases of
``insider trading" or imperfect information and for all other
mixed cases (for instance, we can imagine that we do not observe
the price of the assets but have some extra information which is
not contained in the filtration induced by the processes of
exchange rates. Observe that, if we know that the price of some
asset will double between today and tomorrow, we can  make an
arbitrage without knowing this price - assuming that transaction
costs are reasonable).

In Section \ref{sec abstract formulation}, we  study the
no-arbitrage conditions considered in  \cite{KSR01}, \cite{KSR04}
and \cite{schach} in this abstract setting.   Examples of
application are provided in Section \ref{sec example}

\section{The abstract formulation}\label{sec abstract formulation}

\no Throughout this paper, we fix a finite time horizon $T\in \N$ and
  consider a complete probability space $(\Omega,\Fc,\P)$
supporting  a filtration $\H=(\Hc_t)_{t \in \T}$ with
$\T=\{0,\ldots,T\}$. Importantly, we only assume that
$\Hc_T\subset \Fc$. In particular, most of the processes
considered in this paper need not  be   $\H$-adapted. In all this
paper, inequalities involving random variables must be
understood in the $\Pas$ sense, if it is clear from the context,
 and inclusive relations between
elements of $\Fc$ are assumed to hold up to $\P$-null sets.

\subsection{The model}

\no We consider a closed convex cone $\Ac$ of $\M^d$, $d \ge 1$, and denote by $\F$  the set of continuous maps
 $F$ from $\M^d$
 into $\R^d$ such that\\

\no $\HFa$ :  For $\lambda\ge 0$ and $a \in \M^d$, $\lambda F(a)=F(\lambda a)$ .

\vs3

\no $\HFb$ :  For $\lambda\ge 0$, $\beta\ge 0$  and $a,\;a' \in \Ac$,
$F(\lambda a+\beta a')-\left(\lambda
F(a)+\beta F(a')\right)$
$\in \R^d_+$. \\

\no We then define  $\Fbb$ as the set of $\Fc$-measurable
sequences $F=(F_t)_{t \in \T}$ such that $F_t$ takes a.s. values
in $\F$, for each $t \in \T$. Observe that $\HFa$ implies that
$F(0)=0$.

\vs3

\no   Given $F$ $\in  \Fbb $, we   define  $N(F)$ $=$
$(N_t(F))_{t\in \T}$ and $N^0(F)$ $=$ $(N^0_t(F))_{t\in \T}$ by
    \b*
    N_t(F)\=\left\{F_t(\eta),\;\eta \in L^0(\Ac;\Hc_t)\right\}
    &\And&
    N^0_t(F)\=N_t(F)\cap (-N_t(F))\;.
    \e*
Here, for $E\subset \M^d$ (or $E\subset \R^d$) and a
$\sigma$-algebra $\Gc$ included in $\Fc$, $L^0(E;\Gc)$ denotes
the set of $E$-valued $\Gc$-measurable random variables. For a
process $\xi$ such that $\xi_t\in N_t(F)$ for all $t$ $\in \T$,
we shall simply write $\xi \in N(F)$. We shall similarly write
$\xi \in N^0(F)$ if $\xi_t\in N^0_t(F)$ for all $t$ $\in \T$.

\vs3

\no Given a process $\xi$ with values in $\R^d$, we finally
define
    \b*
    V_t(\xi)\=\sum_{s=0}^t \xi_s &\And&
    A_t(F)\;:=\; \left\{V_t(\xi)-r,\; \xi \in N(F),\; r\in L^0(\R^d_+;\Fc)\right\}\;,
    \;\;\;t\in \T
    \;.
    \e*
Observe that we do not impose that the above processes are
$\H$-adapted:  $F_t$, $\xi_t$ and $V_t(\xi_t)$ need not   be
$\Hc_t$-measurable.

\begin{Remark}\label{rem ex F}{\rm  In financial applications, $F^i_t(\eta_t)$ will correspond to the change in the number of units
of asset $i$ held in the portfolio $V(\xi)$ at time $t$. This
results from the different exchanges $\eta_t^{ij}$ and
$\eta^{ji}_t$ made between the $i$-th asset and the other $j$-th
assets, under the self-financing condition and after paying the transaction costs. 
In this case, $A_t(F)$ stands for the set of contingent
claims, labeled in physical units, that can be super-hedged by
trading up to time $t$ and starting with an initial endowment
equal to $0$. This formalism applies to model \reff{eq modele 1} with $\Ac=\M^d$ and
    \b*
        F^i_t(\eta_t)(\omega)&=& \sum_{j\le d}
    \eta_t^{ji}(1+\lambda^{ij}_t(\omega)\1_{\eta_t^{ij}<0})
    -\eta_t^{ij}\tau^{ij}_t(\omega)(1+\lambda^{ij}_t(\omega)\1_{\eta_t^{ij}>0})\;\;,\;i\le d\;.
    \e*
This model will be further discussed in Section
\ref{sec example}. }
\end{Remark}

\no In this section, we provide sufficient conditions under
which  $A_T(F)$ is   closed in probability and study   abstract
versions of the no-arbitrage conditions of  \cite{KSR01},
\cite{KSR04} and \cite{schach}.

\subsection{Sufficient conditions for the closedness of $A_T(F)$}

\no In all this subsection, we shall assume that the sequence of random maps $F$ satisfies the following conditions: \\

\no $\HKP$: For each $\xi$ and $\tilde \xi$ in $N(F)$,
$V_T(\xi)+V_T(\tilde \xi) \in L^0(\R^d_+;\Fc)$
implies that $\xi \in N^0(F)$ and $V_T(\xi)+V_T(\tilde \xi)=0$.\\

\no $\HNF$: For $t\in \T$ and $\eta \in L^0(\Ac;\Hc_t)$,
$F_t(\eta)\in N^0_t(F) \Rightarrow F_t(-\eta)=-F_t(\eta)$
and $-\eta \in L^0(\Ac;\Hc_t)$. \\

\no We call the first condition $\HKP$ as ``key property"   as it
results from what was called  ``key Lemma" in \cite{KSR04}, see
condition \iii in \cite{KSR04} and Lemma 3 in \cite{KSR01}. In
Subsection \ref{subsec na abstrait}, we shall provide sufficient
conditions   for this property to hold.

\vspace{3mm}

\no In financial models with transaction costs, the second
condition can be understood as follows: $\xi_t:=F_t(\eta_t)\in
N^0_t(F)$ means that the exchange $\xi_t$ is reversible, i.e.
starting with the endowment $\xi_t$ we can make immediate
exchanges so as to come back to $0$. Intuitively, this means
that  $\eta_t$ corresponds to exchanges between assets that can
be exchanged freely, i.e. without paying transaction costs. In
this case, we should be able to do the opposite operation,
$-\eta_t$, to reverse these transactions. In the formalism of
\cite{KSR01} and \cite{KSR04} such an assumption is not required
and the only important property is that if $\xi_t \in N^0_t(F)$
then $-\xi_t \in N_t(F)$, which, in their setting,   implies that
$-\xi_t$ is also an admissible exchange. Since, in our case,
$-\xi_t$ may not be $\Hc_t$-measurable, we need to rewrite it as
some $F_t(\tilde \eta_t)$ for some suitable $\tilde \eta_t \in
L^0(\Ac;\Hc_t)$. In view of the above discussion, it is natural
to assume
that such a $\tilde \eta_t$ should be simply given by $-\eta_t$. \\

\no The aim of this section is to show that it implies the
closedness (in probability) of the set $A_T(F)$.

\vspace{3mm}

\no For the reader's convenience, we  first recall the following
Lemma whose proof can be found in \cite{teachers note}.

\begin{Lemma}\label{lem liminf finie} Set $\Gc \subset \Fc$ and $E\subset \R^{d}$.
Let $(\eta^n)_{n\ge 1}$ be a sequence in $L^0(E;\Gc)$. Set $\tilde
\Omega$ $:=$ $\{\liminf_{n\to \infty}$ $ \|\eta^n\| < \infty\}$.
Then, there is an increasing sequence of random variables
$(\tau(n))_{n\ge 1}$ in $L^0(\N;\Gc)$ such that $\tau(n)\to
\infty$ a.s. and $\eta^{\tau(n)}\1_{\tilde \Omega}$ converges a.s. to   $\eta^*\1_{\tilde \Omega}$
for some  $\eta^*\in L^0(E;\Gc)$.
\end{Lemma}

\no  In the following, we shall denote by  $\L^0(\Ac;\H)$  the
set of $\Ac$-valued  $\H$-adapted processes.

\begin{Proposition}\label{prop AT ferme} Fix $F$ $\in \Fbb$ such that $\HKP$ and $\HNF$ hold. Then,
$A_T(F)$ is closed in probability.
\end{Proposition}

\proof Let us define $A_{t,T}:=\{\sum_{s=t}^T \xi_s-r,\; \xi \in
N(F),\;r\in L^0(\R^d_+;\Fc)\}$, $t\in \T$. We claim that $A_{T,T}$ is closed in
probability (see 3. below) and use an inductive argument. We
assume that $A_{t+1,T}$ is closed in probability for some $t\le
T-1$ and show that $A_{t,T}$ is closed too. Let $(g^n)_{n\ge 1}$
be a sequence in $A_{t,T}$ which converges a.s. to some $g$ $\in$
$L^0(\R^d;\Fc)$. We have to show that $g$ $\in A_{t,T}$. Let
$(\eta^n,r^n)_{n\ge 1}$ be a sequence in $\L^0(\Ac;\H)\x
L^0(\R^d_+;\Fc)$ such that
    \be\label{eq VTn sur rep gn}
    V_T(\xi^n)-r^n=g^n
    \ee
with $\xi^n:=F(\eta^n)$ and $\eta^n=0$ on $\{0,\ldots,t-1\}$. Set
$\alpha^n$ $:=$ $\|\eta^n_t\|$ and $B:=\{\liminf_{n\to \infty}
\alpha^n<\infty\}$. Since $B$ $\in \Hc_t$, we can work separately
on $B$ and $B^c$,  by considering the two sequences
$(\eta^n\1_B,r^n\1_B)_{n\ge 1}$ and $(\eta^n\1_{B^c},r^n\1_{B^c})_{n\ge 1}$,  and
therefore do as if either $\Pro{B}=1$ or $\Pro{B}=0$.

\no {\bf 1.} If $\Pro{B}=1$, then, by Lemma \ref{lem liminf
finie}, there is a random sequence $(\tau(n))_{n\ge 1}$ in
$L^0(\N;\Hc_t)$ such that $\tau(n)\to \infty$ a.s. and
$\eta^{\tau(n)}_t$ converges a.s. to some $\eta^*_t \in
L^0(\Ac;\Hc_t)$. Then, by a.s. continuity of $F_t$,
$F_t(\eta^{\tau(n)}_t)$ converges to $F_t(\eta^*_t)$. Since by
construction $g^{\tau(n)}-F_t(\eta^{\tau(n)}_t)$  $\in
A_{t+1,T}(F)$, and, by assumption, the later is closed in
probability, we can find some $\tilde \xi \in N(F)$ such that
$\tilde \xi=0$ on $\{0,\ldots,t\}$ and $\sum_{s=t+1}^T \tilde
\xi_s=g-F_t(\eta^*_t)$. Since $F_t(\eta^*_t) \in N_t(F)$, this
shows that $g$ $\in A_{t,T}$.

\no {\bf 2.} If $\Pro{B}=0$ then we set $\bar
\eta^n:=\eta^n/(\alpha^n\vee 1)$. Since $\liminf_{n\to \infty}
\|\bar \eta^n_t\|<\infty$ a.s., we can assume (after possibly
passing to a $\Hc_t$-measurable random subsequence as above) that
$\bar \eta^n_t$ converges a.s. to some element of
$L^0(\Ac;\Hc_t)$. Arguing as above, using $\HFa$ and observing
that $g^n/(\alpha^n\vee 1)$ converges a.s. to $0$, we can find
some $\bar \eta \in \L^0(\Ac;\H)$, such that $\bar \eta^n_t\to
\bar \eta_t$, and $\bar r \in L^0(\R^d_+;\Fc)$
 for which
    \be\label{eq sum F etabar=0}
    \sum_{s=t}^T F_s(\bar \eta_s) -\bar r&=&0 \;.
    \ee
From $\HKP$, we deduce that $\bar r=0$ and $\bar \xi_s:=F_s(\bar
\eta_s)$ $\in N^0_s(F)$ for all $s\ge t$. Since $\|\bar \eta_t\|=
1$, there is partition of $B$ into (possibly empty) disjoint sets
$(B_{ij})_{ i,j\le d}$   such that $B_{ij}\subset \{\bar
\eta^{ij}_t\ne 0\}$. We then define $\tilde
\eta^n_s:=\sum_{i,j\le d} (\eta^n_s -\beta_n^{ij}\bar
\eta_s)\1_{B_{ij}}\1_{s\ge t}$ where
$\beta_n^{ij}=(\eta^n_t)^{ij}/\bar \eta^{ij}_t$ on $B_{ij}$ and
$\beta_n^{ij}=0$ on $B_{ij}^c$. Set
$C_{ij}=\{\beta_n^{ij}=|\beta_n^{ij}|\}$ $\cap$ $B_{ij}$ and $\tilde C_{ij}=\{\beta_n^{ij}=-|\beta_n^{ij}|\}$ $\cap$ $B_{ij}$.  By
$\HFa$, $\HFb$, $\HNF$, \reff{eq sum F etabar=0} and the fact
that $\bar r=0$, we get
    \b*
    V_T(F(\tilde \eta^n))
    &=&\sum_{i,j\le d}
    \sum_{s=t}^T  F_s(\eta^n_s)+|\beta_n^{ij}| F_s(\bar \eta_s(\1_{\tilde C_{ij} }-\1_{C_{ij}}))
    +\check r^n 
    \\
    &=&\sum_{i,j\le d}
    \sum_{s=t}^T  F_s(\eta^n_s)
    +|\beta_n^{ij}| \left(F_s(\bar \eta_s)\1_{\tilde C_{ij}  }-F_s(\bar \eta_s)\1_{C_{ij}}\right)
    +\check r^n 
    \\
    &=&
    g^n+r^n+\check r^n   \;,
    \e*
for some sequence  $(\check r^n )_{n\ge 1}$  in $L^0(\R^d_+;\Fc)$. Hence, we
have constructed a new sequence $(\tilde \xi^n:=F(\tilde
\eta^n),\tilde r^n:=r^n+\check r^n)_{n\ge 1}$ for which \reff{eq
VTn sur rep gn} holds and $(\tilde \eta^n_t)^{ij}=0$ on $B_{ij}$.
Repeating this argument recursively on the different $B_{ij}$'s
and arguing as in \cite{KSR01}, we can finally obtain, in a
finite number of operations, a sequence  $(\hat \eta^n)_{n\ge 1}$
in $\L^0(\Ac;\H)$ such that $\liminf_{n\to \infty} \|\hat
\eta^n_t\| <\infty$ a.s. and $\sum_{s=t}^T F_s(\hat
\eta^n_s)=g^n+\hat r^n$, for some sequence $(\hat r^n)_{n\ge 1}$
in $L^0(\R^d_+;\Fc)$. Applying the argument of 1. above then
concludes the proof.

\no {\bf 3.} The fact that $A_{T,T}$ is closed in probability is
obtained by similar   arguments.  Given a sequence $(g^n)_{n\ge
1}$  in $A_{T,T}$ which converges a.s. to some $g$ $\in$
$L^0(\R^d;\Fc)$, we consider a sequence $(\eta^n_T,r^n)_{n\ge 1}$
in $L^0(\Ac;\Hc_T)\x L^0(\R^d_+;\Fc)$ such that
$F_T(\eta^n_T)-r^n=g^n$. Considering separately the event sets
$\{\liminf_{n\to \infty} \|\eta^n_T\|<\infty\}$ and
$\{\liminf_{n\to \infty} \|\eta^n_T\|=\infty\}$ as in 1. and 2.,
we can construct a new sequence $(\hat \eta^n_T,\hat  r^n)_{n\ge
1}$ such that $F_T(\hat  \eta^n_T)-\hat  r^n=g^n$ and
$\liminf_{n\to \infty} \|\hat \eta^n_T\|<\infty$. By possibly
passing to a random subsequence, we can then assume that $\hat
\eta^n_T$ converges a.s. to some $\hat \eta_T$ $\in
L^0(\Ac;\Hc_T)$ and therefore $\hat  r^n$ converges  to some
$\hat r \in L^0(\R^d_+;\Fc)$ for which $F_T(\hat \eta_T)-\hat
r=g$. \ep

\subsection{Abstract weak no-arbitrage property}

\no In this section, we   use Proposition \ref{prop AT ferme} to
provide a dual characterization of the {\sl weak no-arbitrage}
condition studied in \cite{KabStri} and \cite{KSR01}, see also the
references therein,

\vspace{3mm}

\no $\NAw$ : \hspace{1mm} $A_T( F)\cap L^0(\R^d_+;\Fc)=\{0\}$ .

\vspace{3mm}

\no As observed in \cite{schach}, in financial models, it
corresponds to the usual no-arbitrage condition. Here, we keep
the notations of \cite{KabStri} and \cite{KSR01} to enhance the
difference with the notions of {\sl strict} no-arbitrage and {\sl
robust} no-arbitrage that we shall consider in Subsection
\ref{subsec na abstrait}.

\vspace{2mm}

\no We denote by $e_{ij}$ the element of $\M^d_+$ whose component
$(i,j)$ is equal to one and all others are equal to $0$, $i,j\le
d$. In addition to $\HKP$, we make the following assumption on
$\Ac$.

\vs3

\no $\HL$~:~ 1. $F(\delta e_{ij})=0$ if $\delta e_{ij}$ $\notin
\Ac$, $\delta \in \{-1,1\}$, $i,j\le d$.

\hspace{5mm} 2. For $\eta$   in $L^0(\M^d;\Fc)$,  $F(\eta) =
\sum_{i,j\le d} (\eta^{ij})^+ F(e_{ij})+(\eta^{ij})^- F(-e_{ij})$.

\vspace{2mm}

\no Here, $x^+$ and $x^-$ stands for the positive   and negative
parts of $x$.  Condition 1. can be viewed as a convention. The
reason for imposing this assumption will be clear in Section
\ref{sec example}. In the examples of  Section
\ref{sec example}, $e_{ij}$ (resp. $-e_{ij}$) will correspond to a transfer of units of asset 
$i$ so as obtain (resp. get rid of) one unit of asset $j$. Since an order, $\eta$, can be viewed as a composition of single transfers of the form $e_{ij}$ or $-e_{ij}$, condition 2. simply means that the induced changes $F_t(\eta)$ in the portfolio should correspond to the combination of the changes $F_t(e_{ij})$ and $F_t(-e_{ij})$ associated to these single transfers.   

\vs3

\no Observe from $\HFa$, $\HFb$ that, for all $i,j,k\le d$,
    \be\label{eq Fe}
     F^k(e_{ji})\le -F^k(-e_{ji}) \;\mbox{ if } \;(e_{ji},-e_{ji})\in
     \Ac\x\Ac\;,
     \ee
since $F(e_{ji}-e_{ji})=F(0)=0$.

\vs3

\no We shall also assume in the sequel that
    \be\label{eq Feij in L1}
    F_t(e_{ij}) \And F_t(-e_{ij})   \in L^1(\R^d;\Fc) &\pourtout& i,j\le d \;\And\; t \in \T\;.
    \ee
Here, $L^1(\R^d;\Fc)$ denotes the set of $\P$-integrable elements
of $L^0(\R^d;\Fc)$. 

\begin{Remark} {\rm  Observe that we can always reduce to this case by passing to the equivalent
probability measure whose density with respect to $\P$ is defined by $H/ \Esp{H}$
with $H:=\exp(-\sum_{i,j\le d}  \sum_{t\in \T} \|F_t(e_{ij})\|+ \|F_t(-e_{ij})\|)$.
}
\end{Remark}

\begin{Remark}\label{rem eij}{If $F\in \Fbb$ satisfies $\HL$, then it is completely
characterized by the family $\{F(e_{ij}),F(-e_{ij}) \}_{i,j\le d}$.
}
\end{Remark}

\subsubsection{Dual characterization of $\NAw$ under $\HKP$}

\no  For $Z \in  L^\infty(\R^d;\Fc)$, the set of bounded random
variables in $L^0(\R^d;\Fc)$, and $\eta \in \L^0(\Ac;\H)$, we
define
    \b*
    \bar F_t(\eta_t;Z)&:=&\Esp{Z\cdot F_t(\eta_t)~|~\Hc_t} \;\;, \;\;\; t\in \T\;.
    \e*
Here ``$\cdot$" denotes the natural scalar  product of $\R^d$. By
\reff{eq Feij in L1} and $\HL$, these conditional expectations are
well defined.

\vs3

\no We then define $\Dc(F)$ as the set of elements $Z$ of
$L^\infty(\R^d;\Fc)$ satisfying $Z^i>0$   for all $i\le d$ and
such that for all $\eta \in  \L^0(\Ac;\H)$ and $t \in \T$

$\DFa$: $\bar F_t(\eta_t;Z) \le 0$.

$\DFb$:  $F_t(\eta_t)\1_{\bar F_t(\eta_t;Z)=0} \in N^0_t(F)$.

\begin{Theorem}\label{thm separation} Let $F\in \Fbb$ be such that  $\HL$
holds. Then, $\Dc(F) \ne \emptyset$ $\Rightarrow$ $\NAw$.  If
moreover $\HNF$ and $\HKP$ hold, then $\NAw$ $\Rightarrow$
$\Dc(F) \ne \emptyset$.
\end{Theorem}

\no The proof will be provided in the next subsection.

\vs3

\no In order to relate the above result to the literature, we now
provide an alternative characterization of the set $\Dc(F)$. To $Z
\in L^\infty((0,\infty)^d;\Fc)$, we associate the $\H$-martingale
$\bar Z$ defined by $\bar Z_t$ $=$ $\Esp{Z~|~\Hc_t}$, $t\in \T$.
Then, to $F\in \Fbb$ satisfying $\HL$, we associate  $\hat
F(\cdot;Z)$ defined as the element of $\Fbb$ satisfying $\HL$ and
    \b*
    \hat F^k_t(\delta e_{ij};Z)&=&  \Esp{Z^k F^k_t(\delta e_{ij})~|~\Hc_t}/\bar Z^k_t \;\;,\;\;
    i,j,k\le d\;, \;t\in \T\;,\;\delta \in \{-1,1\}\;,
    \e*
see Remark \ref{rem eij}. Observe that $\bar Z \cdot \hat
F(\cdot;Z)$ $=$ $\bar  F(\cdot;Z)$. Given $i,j\le d$, we then introduce the
sequence of random convex cones $\hat K^{ij}(F,Z)$ $=$ $(\hat
K^{ij}_t(F,Z))_{t\in \T}$ defined by
    \b*
    \hat K^{ij}_t(F,Z)(\omega)&=&\mbox{cone}\{-\hat F_t(e_{ij};Z)(\omega)\;,\;-\hat F_t(-e_{ij};Z)(\omega)\}+\R^d_+\;,
    \e*
where, for  $E\subset \R^d$, cone$\{E\}$ is the smallest closed
convex cone that contains $E$. We also define the sequence $\hat
K^{ij*}(F,Z)$ $=$ $(\hat K^{ij*}_t(F,Z))_{t\in \T}$ by
    \b*
    \hat K^{ij*}_t(F,Z)(\omega)&=&\{y\in \R^d~:~x\cdot y\ge 0,\;\pourtout x\in \hat K^{ij}_t(F,Z)(\omega)\}\;.
    \e*
\no In the case of perfect information, i.e. $F$ is $\H$-adapted,
$\hat K(F,Z):=\sum_{i,j\le d} \hat K^{ij}(F,Z)$ coincides with
the sequence of random ``solvency" cones defined in \cite{KSR01}
and  \cite{KSR04}. The following proposition  combined with
Theorem \ref{thm separation} then extends the results of
\cite{KSR01}, \cite{KSR04} and \cite{schach} to our context, see
also   Remark \ref{rem prop carct Dc} below.

\begin{Proposition}\label{prop caract Dc} Let $F\in \Fbb$ be such
that $\HNF$ and $\HL$ hold. Then,  $\Dc(F)$ is the set of   elements
$Z$ of $L^\infty((0,\infty)^d;\Fc)$ such that  $\bar Z_t \in \bigcap_{i,j\le d} \;\ri{\hat K^{ij*}_t(F,Z)}$ a.s.  for all $t\in \T$.
\end{Proposition}

\proof Let $\Dc'(F)$ denote the set of   elements
$Z$ of $L^\infty((0,\infty)^d;\Fc)$ such that  $\bar Z_t \in \bigcap_{i,j\le d} \;\ri{\hat K^{ij*}_t(F,Z)}$,  for all $t\in \T$.

\no {\bf 1.} We fix  $t\in \T$. Since $\bar Z \cdot \hat
F(\cdot;Z)$ $=$ $\bar F(\cdot;Z)$, it follows that $\bar
F_t(\delta e_{ij};Z) \le 0$ for all $\delta \in \{-1,1\}$ is
equivalent to $\bar Z_t \in \hat K^{ij*}_t(F,Z)$, for all $i,j\le
d$.

\no {\bf 2.}  Assume that $\Dc(F)\ne \emptyset$, fix $Z \in
\Dc(F)$ and  $i,j\le d$. Set $B:=\{ \bar Z_t \notin \ri{\hat
K^{ij*}_t(F,Z)} \}$. If $\Pro{B}>0$, we can find some $\hat \xi$
in $L^0(\R^d;\Hc_t)$ with values in $(-\hat K^{ij}_t(F,Z))
\setminus \hat K^{ij}_t(F,Z)$ on $B$ such that $\hat \xi \cdot
\bar Z_t=0$. Since $\hat \xi \in -\hat K^{ij}_t(F,Z)$, there is
some $(\eta,r) \in L^0(\Ac;\Hc_t)$ $\x L^0(\R^d_+;\Fc)$ such that
$\eta^{kl}$ $=$ $0$ if $(k,l) \ne (i,j)$ and $\hat \xi=\hat
F_t(\eta;Z)-r$, recall $\HFb$. By $\DFa$, it satisfies $\bar
F_t(\eta;Z)=0$ and $r=0$. Set $\xi:= F_t(\eta)$. We claim that
$\xi$ $\notin N^0_t(F)$, which, in view of $\DFb$, leads to a
contradiction. To see this observe from $\HNF$ that $\hat
F_t(\eta;Z)$ $=$ $- \hat F_t(-\eta;Z)$ whenever $\xi$ $\in
N^0_t(F)$. By $\HL$, this implies that $\hat \xi \in  \hat K^{ij}_t(F,Z)$
on $B$, a contradiction too. Hence $\Pro{B}=0$. This shows that
$\Dc(F)\subset \Dc'(F)$.

\no {\bf 3.} Assume that $\Dc'(F)\ne \emptyset$ and fix $Z \in
\Dc'(F)$. In view of 1. and $\HL$, it remains to show that if
$\xi$ $\in N_t(F)$ is such that $\Esp{Z\cdot \xi~|~\Hc_t}=0$,
then $\xi \in N^0_t(F)$. Set $\eta \in L^0(\Ac;\Hc_t)$ such that
$ \xi= F_t(\eta)$. Since   $\bar F_t(\eta;Z)=0$, it follows from
$\HL$ that
    \be\label{eq egale 0}
    0&=&\sum_{i,j\le d}  (\eta^{ij})^+   \bar F_t(e_{ij};Z)+
     (\eta^{ij})^-  \bar F_t(-e_{ij};Z)  \;.
    \ee
Since $\bar Z_t \in \bigcap_{i,j\le d} \hat K^{ij*}_t(F,Z) $, we deduce that
    \be\label{eq etaij0 on}
    \eta^{ij} =0 &\mbox{ on } \{\bar F_t(e_{ij};Z)<0\}\cap \{\bar F_t(-e_{ij};Z)<0\}\;.
    \ee
We claim that
    \be\label{eq claim polaire}
    \{\bar F_t(e_{ij};Z)=0\}&=&\{\bar F_t(-e_{ij};Z)=0\} \subset\{F_t(e_{ij})=-F_t(-e_{ij})\}\;,\;i,j\le d\;.
    \ee
In view of \reff{eq etaij0 on} and $\HL$,
this implies that $F_t(\eta)=-F_t(-\eta)$ and therefore $\xi \in N^0_t(F)$.

\no It remains to prove \reff{eq claim polaire}. Fix $i,j\le d$.
Since $\bar Z_t  \in \ri{\hat K^{ij*}_t(F,Z)}$, we must have
$\{\bar F_t(e_{ij};Z)=0\}$ $=$ $\{\bar F_t(-e_{ij};Z)=0\}$ $=:$
$B_{ij}$. If $(e_{ij},-e_{ij})\in \Ac\x \Ac$, then \reff{eq Fe}
implies  that $F_t(e_{ij})+F_t(-e_{ij})=0$ on $B_{ij}$, recall
that $Z$ has a.s. positive components. If $e_{ij}\notin \Ac$,
then $F_t(e_{ij})=0$, recall $\HL$, and $\bar F_t(-e_{ij};Z)=0$
implies $F_t(-e_{ij})=0$ since otherwise $\bar Z_t$ would not
take values in $\ri{\bar K^{ij*}_t(F,Z)}$ a.s. Similarly, if $-e_{ij}\notin
\Ac$, then $F_t(-e_{ij})=0$ and $\bar F_t(e_{ij};Z)=0$ implies
$F_t(e_{ij})=0$.  \ep

\begin{Remark}\label{rem prop carct Dc}{\rm Under the assumptions of Proposition \ref{prop caract Dc},   $\bigcap_{i,j\le d} \;\ri{\hat K^{ij*}_t(F,Z)}$ is a.s. non-empty whenever $\Dc(F)\ne \emptyset$. It follows that
    \b*
    \ri{ \bigcap_{i,j\le d}  \hat K^{ij*}_t(F,Z)}&=&\bigcap_{i,j\le d} \;\ri{\hat K^{ij*}_t(F,Z)}
    \;.
    \e*
Recalling that $\hat K_t(F,Z):=\sum_{i,j\le d} \hat K^{ij}_t(F,Z)$, we then have
        \b*
        \ri{   \hat K^{*}_t(F,Z)}&=&\bigcap_{i,j\le d} \;\ri{\hat K^{ij*}_t(F,Z)}\;,
        \e*
where
        \b*
    \hat K^{*}_t(F,Z)(\omega)&=&\{y\in \R^d~:~x\cdot y\ge 0,\;\pourtout x\in \hat K_t(F,Z)(\omega)\}\;.
    \e*
}
\end{Remark}

\subsubsection{Proof of  Theorem \ref{thm separation}}

\no The two following Lemmas prepare for   the proof of Theorem
\ref{thm separation} which will be concluded at the end of this subsection.

\begin{Lemma}\label{lem pour Dc non vide implique NA} Let $F\in \Fbb$ be such that $\HL$  holds.
Assume that $\Dc(F)\ne \emptyset$.  Then, for all $g \in A_T(F)$ and $Z \in \Dc(F)$ such that $\Esp{Z\cdot g~|~\Hc_T}^- \in L^1(\R;\Fc)$, 
$\Esp{Z\cdot g} \le 0$.
\end{Lemma}

\proof We use a resursive agument as in \cite{KSR01}.  If $g \in A_0(F)$ then $g=F_0(\eta_0)-r$ for some $\eta_0 \in L^0(\Ac;\Hc_0)$ and $r\in L^0(\R^d_+;\Fc)$. By $\DFa$, we have $\Esp{Z\cdot g~|~\Hc_0}\le 
\bar F_0(\eta_0;Z)$ $\le 0$. Next assume that for $g \in A_{t-1}(F)$ such that $\Esp{Z\cdot g~|~\Hc_{t-1}}^-\in L^1(\R;\Fc)$ we have $\Esp{Z\cdot g}\le 0$, for some $0<t\le T$. Then, if $g=\sum_{s=0}^t F_s(\eta_s)-r$ for some $\eta \in L^0(\Ac;\H)$ and $r\in L^0(\R^d_+;\Fc)$, we have $Z\cdot g$ $\le$ $Z\cdot \sum_{s=0}^t F_s(\eta_s)$ and, by $\DFa$, $\Esp{Z\cdot \sum_{s=0}^{t-1} F_s(\eta_s)~|~\Hc_t}\ge  -\Esp{Z\cdot \sum_{s=0}^{t} F_s(\eta_s)~|~\Hc_t}^{-}$. It follows that $\Esp{Z\cdot \sum_{s=0}^{t-1} F_s(\eta_s)~|~\Hc_{t-1}}^- \in L^1(\R;\Fc)$ and therefore $\Esp{Z\cdot \sum_{s=0}^{t-1} F_s(\eta_s)}\le 0$. Since by $\DFa$, $\bar F_t(\eta_t;Z)\le 0$, it follows that $\Esp{Z\cdot g}\le 0$.  Observe that we
have no problem in defining the above conditional expectations
thanks to \reff{eq Feij in L1} and $\HL$.
\ep

\begin{Lemma}\label{lem Zzeta} Let $F\in \Fbb$ be such that $\NAw$, $\HKP$, $\HNF$ and $\HL$ hold. Then, for all
$t \in \T$ and $\mu \in L^0(\Ac;\Hc_t)$, there is $Z^\mu \in
L^\infty(\R^d;\Fc)$ with $(Z^\mu)^i>0$   for all $i\le d$ such
that

\i $\bar F_s(\eta_s;Z^\mu) \le 0$   for all $\eta \in
\L^0(\Ac;\H)$ and $s \in \T$

\ii $F_t(\mu)\1_{\bar F_t(\mu;Z^\mu)=0} \in N^0_t(F)$.
\end{Lemma}

\proof We follow the argument of Lemma 4 in \cite{KSR04}. Observe
from $\HFa$ and $\HFb$ that $A^1_T(F):= A_T(F)\cap L^1(\R^d;\Fc)$
is a convex cone which contains $-L^1(\R^d_+;\Fc)$. Since it is closed in $L^1(\R^d;\Fc)$, see
Proposition \ref{prop AT ferme}, and satisfies $A^1_T(F)\cap
L^1(\R^d_+;\Fc)$ $=$ $\{0\}$, see $\NAw$, we deduce from the
Hahn-Banach separation theorem together with a classical
exhaustion argument, see e.g. Section 3 in \cite{Sinfinite}, that
there is some $Z \in L^\infty(\R^d;\Fc)$ with $Z^i>0$ for all
$i\le d$ such that $\Esp{Z\cdot g}\le 0$ for all $g$ $\in
A^1_T(F)$. Let $\Zc$ denote the set of such random variables $Z$.

\no {\bf 1.} It is clear that \i holds for all $Z$ $\in  \Zc$.
Indeed, assume that for some $\eta  \in  \L^0(\Ac;\H)$ and $s \in
\T$, $B:=\{\bar F_s(\eta_s;Z)> 0\}$ has positive probability.  Set
$\tilde g:=H_sF_s(\eta_s)\1_B$ with $H_s:=\exp(-\|\eta_s\|) \in
L^0((0,\infty);\Hc_s)$.  By $\HFa$, $H_s F_s(\eta_s)\1_B$ $=$
$F_s(H_s \eta_s\1_B)$ so that, by \reff{eq Feij in L1} and $\HL$,
$\tilde g \in A^1_T(F)$. Since $\Esp{ Z\cdot \tilde g}> 0$,   we
get a contradiction to the definition of $\Zc$.

\no {\bf 2.} By the same argument as in Lemma 4 in \cite{KSR04},
we can find some $Z^\mu$ such that $\Pro{\bar F_t(\mu;Z^\mu)<0}$
$=$ $\max_{Z \in \Zc} \Pro{\bar F_t(\mu;Z)<0}$. Set $B:=\{\bar
F_t(\mu;Z^\mu)=0 \}$ and $B_k:=B \cap \{ \|\mu\|\le k\}$, $k\in
\N$. We claim that if \ii fails for $(\mu,Z^\mu)$ then
$-F_t(\mu\1_{B_k}) \notin$ $A^1_T(F)$ for some $k>0$. Indeed,
otherwise, for all $k>0$,  we could find some $ \eta_k$ $\in
\L^0(\Ac;\H)$ and $r_k \in L^0(\R^d_+;\Fc)$ such that $V_T(F(
\eta_k))=-F_t(\mu\1_{B_k})+r_k$ $ \in  A^1_T(F)$, so that $V_T(F(
\eta_k))+F_t(\mu\1_{B_k})$  $\in L^0(\R^d_+;\Fc)$. By  $\HKP$,
this would imply that $F_t( \mu\1_{B_k}) \in N^0_t(F)$, so that,
by $\HNF$,  $F_t( \mu\1_{B_k})=-F_t( -\mu\1_{B_k})$. Sending $k
\to \infty$, we would then get $F_t( \mu\1_{B})=-F_t(
-\mu\1_{B})$, showing that $F_t( \mu\1_{B}) \in -N_t(F)$, a
contradiction. Hence, if \ii fails  $-F_t(\mu\1_{B_k}) \notin$
$A^1_T(F)$ for some $k>0$. Repeating the argument of 1., we can
then find some $Z \in L^\infty(\R^d_+;\Fc)$  such that $\Esp{Z
\cdot g}\le 0  < \Esp{Z \cdot (-F_t(\mu\1_{B_k}) )}$ for all $g
\in A^1_T(F)$. Taking $\tilde Z=Z+Z^\mu$, we obtain $\Pro{\bar
F_t(\mu;\tilde Z)<0}> \Pro{\bar F_t(\mu; Z^\mu)<0}$, a
contradiction to the definition of $Z^\mu$. This shows that \ii
must hold.
 \ep\\

\no {\bf Proof of Theorem \ref{thm separation}.} {\bf 1.} The
first implication follows from Lemma \ref{lem pour Dc non vide
implique NA} since the elements of $\Dc(F)$ have a.s. positive
entries.

\no {\bf 2.} We now prove the converse implication. Let
$Z^{ij,t}_+$ (resp. $Z^{ij,t}_-$) be an element of
$L^\infty((0,\infty)^d;\Fc)$ such that \i and \ii of Lemma
\ref{lem Zzeta} hold for the process $(e_{ij}\1_{s=t})_{s\in \T}$
(resp. $(-e_{ij}\1_{s=t})_{s\in \T}$), $i,j\le d$ and $t\in \T$.
We claim that
    $
        \hat Z :=  \sum_{t\in \T} \sum_{i,j\le d} Z^{ij,t}_+ + Z^{ij,t}_-\;
    $
belongs to  $\Dc(F)$. Clearly, it satisfies $\DFa$. Fix  $\eta
\in \L^0(\Ac;\Hc_t)$ for some $t\in\T$, and recall from $\HL$ that
    \be\label{eq F=somme Fe}
    F_t(\eta) &=& \sum_{i,j\le d}  (\eta^{ij})^+ F_t(e_{ij}) + (\eta^{ij})^-F_t(-e_{ij})\;.
    \ee
Set $B:=\{\bar F_t(\eta;\hat Z)=0\}$. From  the definition of $(Z^{ij,t}_+,Z^{ij,t}_-)_{i,j,t}$, we deduce that  $(\eta^{ij})^+
F_t(e_{ij})\1_B$ and $(\eta^{ij})^- F_t(-e_{ij})\1_B$ belongs to
$N^0_t(F)$ for all $i,j\le d$. By $\HNF$, $\HL$ and \reff{eq
F=somme Fe}, we then deduce that
    \b*
    -F_t(\eta\1_B)
    &=& \sum_{i,j\le d}  -(\eta^{ij})^+ F_t(e_{ij}) - (\eta^{ij})^- F_t(-e_{ij})\\
    &=& \sum_{i,j\le d}   (\eta^{ij})^+ F_t(-e_{ij}) + (\eta^{ij})^- F_t(e_{ij})\\
    &=&     F_t(-\eta\1_B)
    \e*
so that $F_t(\eta\1_B)=-F_t(-\eta\1_B) \in -N_t(F)$ and therefore
$F_t(\eta\1_B)\in  N^0_t(F)$. Hence, $\hat Z$ satisfies $\DFb$. \ep

\subsection{Strict and robust  no-arbitrage conditions}\label{subsec na abstrait}

\no In this section, we study the other no-arbitrage conditions
considered in \cite{KSR01}, \cite{KSR04} and  \cite{schach}.

\vs2

\no Following \cite{KSR01}, we say that
 $F$ $\in \Fbb$ satisfies the {\sl   strict no-arbitrage condition} if one has
    \b*
  \NAs\;: &&    A_t(F)\cap (-N_t(F)+L^0(\R^d_+;\Fc)) \subset N^0_t(F) \;\;\;\;\pourtout \;t\in \T\;,
    \e*
and that the model has ``efficient frictions" if
    \b*
   \EF  \;:&& N^0_t(F)=\{0\}\;\;\; \pourtout\; t \in \T \;.
    \e*

\no As in \cite{schach}, we also define a {\sl robust} version of
the no-arbitrage property. We  say that  $F$ $\in \Fbb$ satisfies
the {\sl robust no-arbitrage condition},  $\NAr$, if there is some
sequence $ G \in \Fbb$ such   that for all $\eta \in
\L^0(\Ac;\H)$, $t$ $\in \T$ and $i\le d$:
    \b*
  &1.&   G^i_t(\eta_t) \ge F^i_t(\eta_t)
    \\
    &2.& F_t(\eta_t) \notin N^0_t(F) \;\Rightarrow\;  \{\exists \;k\le d \mbox{ such that }  G^k_t(\eta_t)>
    F^k_t(\eta_t)\} \ne \emptyset
 \;\;\;\;
    \\
    &3.&  \mbox{$\NAw$ holds for $G$. }
    \e*

\no In financial models, the last condition can be interpreted as
the existence of a model with slightly lower transaction costs
(for those that are not already equal to $0$) in which  the  {\sl   weak
no-arbitrage condition} still holds, see \cite{schach}.

 \vs3

\no In this section, we first show that these properties imply the
condition $\HKP$ used above. We will then be
able to use Theorem \ref{thm separation} to provide a dual
characterization of the absence of arbitrage opportunities in the
spirit of \cite{KSR01},  \cite{KSR04} and \cite{schach}, see
Theorem \ref{thm caracterisation duale}  below.

\begin{Lemma}\label{lem key lemma}  Let $F\in \Fbb$ be such that  one of the above conditions holds:

\i $\NAr$

\ii $\NAs$ and $\EF$.

\no Then, $\HKP$ holds.
\end{Lemma}

\proof  Set $\xi$ and $\tilde \xi$ in $N(F)$ such that
$V_T(\xi)+V_T(\tilde \xi)$ $ \in L^0(\R^d_+;\Fc)$. Let $\eta$ and
$\tilde \eta$ be  elements of $\L^0(\Ac;\H)$ such that
$\xi=F(\eta)$ and  $\tilde \xi=F(\tilde \eta)$,  set $\bar
\eta:=\eta+\tilde \eta$ and $\bar \xi$ $:=$ $F(\bar \eta)$.

\no {\bf 1.} We start with $\NAr$. Let $ G$ be as in the
definition of $\NAr$  and define   $\bar \xi':=G( \bar \eta)$. By
1. and 2.  of $\NAr$, if for some $t \in \T$    $F_t(\eta) \notin
N^0_t(F)$, then we can find $i\le d$ and $B\in \Fc$ with positive
measure such that  $V^i_T( G(\bar \eta))> V^i_T(\bar \xi)$ on
$B$. By 1. of $\NAr$ and $\HFb$, we then have   $V_T( \bar
\xi')-V_T( \xi)- V_T(\tilde \xi) \in
L^0(\R^d_+;\Fc)\setminus\{0\}$. Since $V_T(\xi)+V_T(\tilde \xi)$
$ \in L^0(\R^d_+;\Fc)$, this leads to a contradiction to  the
fact that $\NAw$ holds for $G$. Hence, $F(\eta) \in N^0(F)$ and
we must have  $V_T( \bar \xi')=0$ so that $V_T( \xi)+ V_T(\tilde
\xi) $ $=$ $0$.

\no {\bf 2.} We now assume that $\NAs$ and $\EF$ hold. Assume
that, for some $t\in \T$, $\xi_t \notin N^0_t(F)$  or $\tilde
\xi_t \notin N^0_t(F)$ and set $t^*$ $:=$ $\max\{t\in \T~:~$
$\xi_t \notin N^0_t(F)$  or $\tilde \xi_t \notin N^0_t(F)   \}$.
Then, by $\EF$ and $\HFb$, $V_{t^*-1}(\bar \xi)$ $=$
$\sum_{s=0}^{t^*-1}\xi_s+ \tilde \xi_s+r$ $=$ $-\xi_{t^*}-\tilde
\xi_{t^*}+r$ $=$ $- \bar \xi_{t^*} +r+r'$ for some $r$, $r'$ in
$L^0(\R^d_+;\Fc)$. This shows that $V_{t^*-1}(\bar \xi) \in
(-N_{t^*}(F)+L^0(\R^d_+;\Fc))\cap A_{t^*}(F)$. By $\NAs$ and
$\EF$, we must have $\bar \xi_{t^*} \in N^0_{t^*}(F)=\{0\}$ and
$r$ $=$ $r'$ $=0$. Hence,
 $\xi_{t^*}=-\tilde \xi_{t^*} \in N^0_{t^*}(F)$,
 thus providing a contradiction to the definition of $t^*$.
 \ep\\

\no Observe that $\NAr$ implies $\NAw$ and that $\NAs$ also
implies $\NAw$ whenever $N^0_T$ $=$ $\{0\}$. In view of Lemma
\ref{lem key lemma}, we can then apply Proposition \ref{prop AT
ferme} and Theorem \ref{thm separation} to deduce that, under
$\HNF$ and $\HL$,  $\NAr$ and ($\NAs$ and $\EF$) both imply that
$A_T(F)$ is closed in probability and that $\Dc(F)$ is non-empty.
Conversely, if $\Dc(F)$ $\ne$ $\emptyset$, on can show that
$\NAs$ and $\NAr$ hold.

\begin{Theorem}\label{thm caracterisation duale} Let $F\in \Fbb$ be such that $\HNF$ and $\HL$
hold.  Then,

\no \i  If either $\NAr$ or {\rm ($\NAs$ and $\EF$)} hold, then
$\Dc(F) \ne \emptyset$ and $A_T(F)$ is closed in probability.

\no \ii If $\Dc(F) \ne \emptyset$ then  $\NAs$ and $\NAr$ hold.
\end{Theorem}

\proof {\bf 1.} Since $\NAr$ implies $\NAw$ and   $\NAs$ also
implies $\NAw$ whenever $\EF$ holds, combining Lemma \ref{lem key
lemma} with Proposition \ref{prop AT ferme} and Theorem \ref{thm
separation} leads to \i.  To show that $\NAs$ holds under $\Dc(F)
\ne \emptyset$, we set $V_t \in A_t(F)$ such that $V_t=$
$-F_t(\tilde \eta)+r$  for some $\tilde \eta \in L^0(\Ac;\Hc_t)$
and $r\in L^0(\R^d_+;\Fc)$. By $\DFa$, $\Esp{Z\cdot V_t~|~\Hc_t}$ $\ge $ $-\bar F_t(\tilde \eta;Z)$ $\ge 0$,  and therefore, by Lemma \ref{lem pour Dc non vide implique NA}, we must have $\Esp{Z\cdot V_t~|~\Hc_t}=0$, $r=0$ and $\bar
F_t(\tilde\eta;Z)=0$ for all $Z\in \Dc(F)$. Then $\DFb$ implies
that $F_t(\tilde \eta)\in N^0_t(F)$.

\no {\bf 2.} We now prove  that $\Dc(F)\ne \emptyset$ implies
$\NAr$. To avoid unnecessary complications, we first consider the
case where $(e_{ji},-e_{ji}) \in \Ac\x\Ac$ for all $i,j\le d$. We
shall explain in 2.d. how to adapt our arguments to the general
case.

\no Fix $Z\in \Dc(F)$ and consider the random variables
    \b*
    \delta_{ji,t}^+ \;:=\; -\bar F_t(e_{ji};Z)
    &\And&
    \delta_{ji,t}^- \;:=\; -\bar F_t(-e_{ji};Z) \;\;,\;\;i,j\le d\;,\;t\in \T\;.
    \e*
It follows from   $\DFa$  that
    \be\label{eq diff delta ok}
      \delta_{ji,t}^+  \ge  0  &\And& \delta_{ji,t}^-\ge 0\;\;,\;\; i,j\le
      d\;,\;t\in \T\;.
    \ee
We claim that, for all $i,j\le d$ and $t\in \T$,
    \be
      && \delta_{ji,t}^+  >  0 \; \And  \; \delta_{ji,t}^-> 0  \;\mbox{
      on }    \{\bar F_t(e_{ji};Z)<0\}=\{\bar F_t(-e_{ji};Z)<0\}\;.
      \label{eq diff > 0 sur bar F<0}
    \ee
Indeed, by construction, we have $\delta_{ji,t}^+  >  0$ on
$\{\bar F_t(e_{ji};Z)$ $<0\}$ and $ \delta_{ji,t}^-> 0$ on $\{\bar
F_t(-e_{ji};Z)$ $<0\}$. Now,  set $B_+:=\{\bar F_t(e_{ji};Z)=0\}$ and
$B_-:=\{\bar F_t(-e_{ji};Z)=0\}$. From $\DFb$ and $\HNF$, we
deduce that $F_t(e_{ji}\1_{B_+})$ $=$ $-F_t(-e_{ji}\1_{B_+})$ so that $\bar F_t(-e_{ji};Z)$ $=$ $0$ on
$B_+$. This shows that $B_+\subset B_-$. Similarly, we can show
the converse inclusion, which implies \reff{eq diff > 0 sur bar
F<0}.

\no We can now construct $G$. For all $i,j,k\le d$, we
set
    \be
    G^k(e_{ji})&=& \left( F^k(e_{ji}) + \delta_{ji,t}^+ /(d\;\bar Z^k_t) \right)\wedge (-F^k(-e_{ji}))  \nonumber
    \\
    G^k(-e_{ji})&=& \left( F^k(-e_{ji}) + \delta_{ji,t}^- /(d\;\bar Z^k_t)   \right)\wedge (-G^k(e_{ji}))
    \;.
    \label{eq def f'}
    \ee
For $x \in \M^d$, we then set
    \b*
     G(x) &=& \sum_{i,j\le d} (x^{ji})^+ G(e_{ji}) +  (x^{ji})^- G(-e_{ji}) \;.
    \e*
It satisfies   $\HFa$. By \reff{eq Fe}, it also satisfies the
condition 1. of $\NAr$, recall \reff{eq diff delta ok}.  It
remains to check that $\HFb$, 2. and 3. of $\NAr$ hold.

\no {\bf 2.a.} We first check $\HFb$. We fix $i,j,k\le d$,
$\alpha\ge\beta\ge 0$. Then, $G^k(\alpha e_{ji} -\beta  e_{ji} )$
$=$   $(\alpha-\beta) G^k( e_{ji})$. By \reff{eq def f'}, it
follows that $G^k(\alpha e_{ji} -\beta  e_{ji} )$ $\ge$   $
\alpha G^k( e_{ji})$ $+ \beta  G^k( -e_{ji})$.  In the case where
$\beta\ge\alpha\ge 0$, we obtain the same result. Since $G$
satisfies $\HL$, this shows that it also satisfies $\HFb$.

\no {\bf 2.b.} We now check 2. of $\NAr$. Set $\eta \in
L^0(\Ac;\Hc_t)$ and $t \in \T$ such that $F_t(\eta)\notin
N^0_t(F)$. We must show that, with positive probability,  we can
find $k\le d$ such that $G^k_t(\eta)>F^k_t(\eta)$. First observe
that we cannot have $\{\eta^{ji}\ne 0\}$ $ \subset$
$\{F_t(e_{ji})=-F_t(-e_{ji})\}$ for all $i,j\le d$ since this
would imply that $F_t(\eta)\in N^0_t(F)$. Hence, there is $(i,j)$
and $k\le d$ such that $B:=\{\eta^{ji}\ne 0\} \cap
\{F^k_t(e_{ji})<-F^k_t(-e_{ji})\}$ $\ne \emptyset$, recall
\reff{eq Fe}. Since, by $\DFa$ and \reff{eq Fe},
$\{F^k_t(e_{ji})<-F^k_t(-e_{ji})\}$ $\subset$ $\{\bar
F_t(e_{ji};Z)<0\}$ $\cup$  $\{\bar F_t(-e_{ji};Z)<0\}$, \reff{eq
diff > 0 sur bar F<0} and \reff{eq def f'} imply that
$G^k_t(\eta)>F^k_t(\eta)$ on $B$.

\no {\bf 2.c.} To check 3. of $\NAr$, it suffices to observe that,
for $\eta \in \L^0(\Ac;\H)$ and $Z \in \Dc(F)$,  we have $\bar
G(\eta;Z)\le 0$. Since $Z$ has a.s. positive entries, the same
arguments as in Lemma \ref{lem pour Dc non vide implique NA}
imply $\NAw$ for $G$.

\no {\bf 2.d.} We now explain how to consider the case where some
$e_{ji}$ or $-e_{ji}$  do not belong to $\Ac$. We assume that,
for some $(i,j)$,
 $e_{ji}$ or $-e_{ji} \in \Ac$, otherwise there is nothing to prove.
We keep the definition of $G$ as above except that in the right
hand-sides of \reff{eq def f'}, we replace  $-F^k(-e_{lm})$ by
$+\infty$ if $-e_{lm} \notin \Ac$ and $-G^k(e_{lm})$ by $+\infty$
if $e_{lm} \notin \Ac$.    Using the convention 1. of $\HL$, we
see that $G$ satisfies $\HFa$ and 1. of $\NAr$. The arguments of
2.a. and 2.c. still hold, so that it also satisfies $\HFb$ and 3.
of $\NAr$. To obtain 2. of $\NAr$, we just recall that
$F^k_t(e_{lm})<0$, for some $k\le d$, implies $\bar
F_t(e_{lm};Z)<0$ whenever $ -e_{lm}\notin \Ac$, see $\DFb$,
$\HNF$ and recall 1. of $\HL$. With this in mind, adapting the
arguments of 2.b. is straightforward.
  \ep

\section{Applications to financial markets with proportional transaction costs}\label{sec example}

\no In this section, we apply the above results to three examples
of discrete time financial markets with proportional transaction
costs.  The first one corresponds to a ``security market" where it
is possible to make transactions only between a ``non-risky
asset" and some ``risky" ones, direct transactions between the
``risky assets" being prohibited. The two other ones correspond to
``currency markets"  where   transactions between all assets
(interpreted as currencies) are possible.  The information of the
financial agent is modeled by the filtration $\H$ and a strategy
is a process $\eta \in \L^0(\Ac;\H)$.

\subsection{Security market}

\no We take the first asset as a num\'eraire and consider an 
$\M^{d}_+$-valued process $\pi$ such that $\pi^{1i}\ge
\pi^{i1}>0$    for all $i,j\le d$ and $\pi^{ii}=1$ for all $i\le
d$. Here, $\pi^{i1}$ must be interpreted as the number of
physical units of asset $1$ one receives when selling one unit of
$i$, and $\pi^{1i}$ as the number of units of asset $1$ one pays
to buy one unit of $i$.  The condition $\pi^{1i}_t\ge \pi^{i1}_t$
is natural since otherwise their would be trivial arbitrages. The
case $\pi^{1i}_t= \pi^{i1}_t$ (resp. $\pi^{1i}_t> \pi^{i1}_t$)
corresponds to the situation with no-friction (resp. with
frictions) between the assets $i$ and $1$.

\no We construct the sequence of random maps $F$ as follows. To
$\rho \in \M^{d}_+$ such that $\rho^{1i}\ge \rho^{i1} >0$, we
associate the map $f(\cdot;\rho)$ from $\M^d$ into $\R^d$ defined
by
    \b*
    f^1(a;\rho)
    \=
    \sum_{i\le d}
    a^{1i}\left(\rho^{i1}\1_{a^{1i}>0} +  \rho^{1i} \1_{a^{1i}<0}\right)
   &\And&
    f^i(a;\rho)
    \=
    -a^{1i} \;\;\;\mbox{ for } i>1
    \;.
    \e*
Then, we set $F_t(\cdot)=f(\cdot;\pi_t)$ for $t\in \T$. For the sake
of simplicity, we take $\Ac=\M^d$. Observe that $\HL$ and $\HFa$
trivially holds, and  that the condition $\pi^{1i}\ge \pi^{i1}$,
$i\le d$, implies $\HFb$.

\vs3

\no If positive, the quantity $\eta^{1i}_t$ corresponds to the
number of units of asset $i$ which are sold in exchange of
$\eta^{1i}_t\pi^{i1}_t$ units of asset $1$. Otherwise
$|\eta^{1i}_t|$ corresponds to the number of units of asset $i$
which are obtained by converting  $|\eta^{1i}_t\pi^{1i}_t|$ units
of asset $1$. The other components of $\eta$ play no role in this
model.

\vs3

\no In order to apply the result of the previous section, we
first check that  $\HNF$  holds  in this model.

\begin{Lemma}  Let $F$ be defined as above, then $\HNF$ holds.
\end{Lemma}

\proof Fix $t$ $\in$ $\T$ and $\eta \in L^0(\M^d;\Hc_t)$  such
that $F_t(\eta)\in N^0_t(F)$. We have to show that
$F_t(-\eta)=-F_t(\eta)$.  By definition, there is $\tilde \eta$
$\in L^0(\M^d;\Hc_t)$ such that $F_t(\eta)$ $=$ $-F_t(\tilde
\eta)$.  Define $S\in L^0((0,\infty)^d;\Fc)$ by $S^i$ $=$
$(\pi_t^{1i}+\pi^{i1}_t)/2$, $i\le d$. Recalling that
$\pi^{11}=1$, direct computation shows that
    \b*
    0\=S\cdot (F_t(\eta )+F_t(\tilde \eta ))
    &=&
     \sum_{i=1}^d
     \eta^{1i}\left(-(\pi_t^{1i}+\pi^{i1}_t)/2+ \pi^{i1}_t \1_{\eta^{1i}>0}+ \pi^{1i}_t \1_{\eta^{1i}<0}\right)
     \\
     &+&
     \sum_{i=1}^d
     \tilde \eta^{1i}\left(-(\pi_t^{1i}+\pi^{i1}_t)/2+ \pi^{i1}_t \1_{\tilde \eta^{1i}>0}+ \pi^{1i}_t \1_{\tilde \eta^{1i}<0}\right)
     \\
     &=&
     \sum_{i=1}^d
     \left(|\eta^{1i}|+|\tilde \eta^{1i}|\right)\left(  \pi_t^{i1}-\pi^{1i}_t\right)/2 \;.
     \e*
Since $\pi^{1i}_t\ge \pi^{i1}_t$ for all $i,j\le d$, this shows that $\eta^{1i}$ is equal to $0$ on
$\{\pi^{1i}_t-\pi_t^{i1}>0\}$ and therefore $F_t(-\eta)=-F_t(\eta)$.
\ep
\\

\no Then, it follows from Theorem \ref{thm caracterisation duale}
that $\NAr$ $\Leftrightarrow$ $\Dc(F)\ne \emptyset$ $\Rightarrow$
$\NAs$ and that the last implication is an equivalence if $\EF$
holds.  We then assume that  $\NAr$ or ($\NAs$ and  $\EF$) hold,
fix $Z \in \Dc(F)$, and define    the process  $\bar \pi$   by
 \b*
 \bar \pi^{ij}_t :=\Esp{Z^1\pi^{ij}_t ~|~\Hc_t}/ \bar Z^1_t\;\;,\;\;i,j\le d\;.
 \e*

\no With this notation, one easily checks that $\bar Z_t \in
\ri{\hat K^{1i*}_t(Z,F)}$   if and only if
    \b*
    \bar Z^1_t \bar \pi^{i1}_t \le \bar Z^i_t \le \bar Z^1_t \bar \pi^{1i}_t\;\;,
    \e*
with strict inequalities   on $\{\bar \pi^{1i}_t>\bar \pi^{i1}_t\}$.

\vs3

\no  Let $\Q$ be the equivalent probability measure defined by
$d\Q/d\P=Z^1/\Esp{Z^1}$. Then, $\bar \pi$ is the optional
projection under $\Q$ of $\pi$ on $\H$, i.e. $\bar \pi_t$ $=$
$\E^{\Q}[\pi_t~|~\Hc_t]$, and there is a $(\Q,\H)$-martingale
$\bar Z/\bar Z^1$ such that each component $i$ evolves in the relative interior of
 the ``estimated" bid-ask spread $[\bar \pi^{i1}_t, \bar \pi^{1i}_t]$. This extends
the discrete-time version of the result of \cite{JKequi}.

\no In the ``no frictions" case, i.e. $\pi^{i1}=\pi^{1i}$, then
$\bar Z^1_t \bar \pi^{i1}_t $ $=$ $\bar Z^i_t$ $=$ $\bar Z^1_t
\bar \pi^{1i}_t$ and we deduce that  there is an equivalent
probability measure under which the optional projection     $\bar
\pi$  of the discounted price processes  $\pi$  on $\H$ are
$(\Q,\H)$-martingales. This is the result of \cite{KSdelayed}.

\subsection{Currency market $\# 1$}\label{subsec currency}

\no We now consider a $(0,\infty)^d$-valued process $S$ which
models the price of the different currencies, before transaction
costs. Then $\tau^{ji}_t=S^i_t/S^j_t$   is the number of units of
asset $j$ that one can exchange at time $t$ against one unit of
asset $i$, before to pay the transaction costs. Transaction costs
are modeled by a process $\lambda$ with values in $\M^d_+$, i.e.
$\lambda^{ji}_t$ is the proportional costs to pay in units of $j$
for an exchange at time $t$ between $j$ and $i$.

\no To construct the sequence of random maps $F$, we first define
the maps $f(\cdot;\rho,\ell)$ from $\M^d$ into $\R^d$  by
    \b*
    f^i(a;\rho,\ell)
    \=
    \sum_{j=1}^d
    a^{ji}\left(1+\ell^{ij}\1_{a^{ji}<0}\right)
    -
    a^{ij}\rho^{ij}\left(1+\ell^{ij}\1_{a^{ij}\ge 0}\right)
    \;,\; \rho,\; \ell \in \M^d_+\;.
    \e*
Then, we set $F_t(\cdot)=f(\cdot;\tau_t,\lambda_t)$, $t\in \T$. For
$\Ac=\M^d$, this corresponds to the model \reff{eq modele 1}
described in the introduction.  Clearly $\HL$, $\HFa$ and $\HFb$
hold.

\no The quantity $\eta^{ij}_t$ corresponds to number of units of asset $j$ which
are obtained by converting     units of asset $i$.
 For such an exchange, the transaction costs are paid in units of asset $i$.\\

\no Here again, we need to check that  $\HNF$   holds  in this
model. For sake of simplicity, we assume that  $\{\lambda^{ij}_t>0\}=\{\lambda^{ji}_t>0\}$ for all $i,j\le d$ and $t\in \T$.

\begin{Lemma}  Let $F$ be defined as above, then $\HNF$  holds.
\end{Lemma}

\proof     Fix $t$ $\in$ $\T$ and $\eta \in L^0(\M^d;\Hc_t)$ such
that $F_t(\eta)\in N^0_t(F)$. We have to show that
$F_t(-\eta)=-F_t(\eta)$.  By definition, there is $\tilde \eta$
$\in L^0(\M^d;\Hc_t)$ such that $F_t(\eta)$ $=$ $-F_t(\tilde
\eta)$. Since $F_t( \eta)+F_t(\tilde \eta)=0$, direct computation
shows that
    \b*
    0\=S_t\cdot (F_t( \eta )+F_t(\tilde \eta ))
    &=&
    -
    \sum_{i,j=1}^d
    | \eta^{ij}| S^j_t \left(\lambda^{ji}_t \1_{ \eta^{ij} < 0}+ \lambda^{ij}_t \1_{ \eta^{ij} >  0}\right)
    \\
    && -
    \sum_{i,j=1}^d
    |\tilde \eta^{ij}| S^j_t \left(\lambda^{ji}_t \1_{\tilde \eta^{ij}< 0}+\lambda^{ij}_t\1_{\tilde \eta^{ij} >  0}\right)
    \;.
    \e*
This shows that $(\eta^{ij})^+=( \eta^{ij})^-=0$   on $\{\lambda^{ij}_t>0\}=\{\lambda^{ji}_t>0\}$
 and   therefore
$F_t(-\eta)=-F_t(\eta)$. \ep
\\

\no It then follows from Theorem \ref{thm caracterisation duale}
that $\NAr$ $\Leftrightarrow$ $\Dc(F)\ne \emptyset$ $\Rightarrow$
$\NAs$ and that the last implication is an equivalence if $\EF$
holds.  We then assume that  $\NAr$ or ($\NAs$ and  $\EF$) hold,
fix $Z \in \Dc(F)$,  and define    the processes  $\bar \tau$ and $\bar \lambda$
by
 \b*
 \bar \tau^{ij}_t :=\Esp{Z^i\tau^{ij}_t\left(1+\lambda^{ij}_t \right) ~|~\Hc_t}/ (\bar Z^i_t (1+ \bar \lambda^{ij}_t))
     &\mbox{ and } & \bar  \lambda^{ij}_t :=  \Esp{Z^i \lambda^{ij}_t  ~|~\Hc_t}/\bar Z^i_t \;.
 \e*

\no With this notation, one easily checks that $\bar Z_t \in
\ri{\hat K^{ij *}_t(Z,F)}$  if and only if
    \b*
    \bar Z^j_t \bar \tau^{ji}_t/(1+\bar \lambda^{ij}_t)
    \le \bar Z^i_t \le \bar Z^j_t \bar \tau^{ji}_t (1+\bar \lambda^{ji}_t)
    \;  \;,
    \e*
with strict inequalities on $\{ \bar \tau^{ji}_t (1+\bar \lambda^{ji}_t)>\bar \tau^{ji}_t/(1+\bar \lambda^{ij}_t)\}$.


\subsection{Currency market $\# 2$}\label{subsec currency 2}

\no The model \reff{eq modele usuel} discussed in the
introduction corresponds to the one presented in the previous
subsection with $f$ defined by
    \b*
    f^i(a;\rho,\ell)
    \=
    \sum_{j=1}^d
    a^{ji} \1_{a^{ji}>0}
    -
    a^{ij}\rho^{ij}\left(1+\ell^{ij}\right)\1_{a^{ij}> 0}
    \;,\; \rho,\; \ell \in \M^d_+\;,
    \e*
i.e. $F$ is defined by $F_t(\cdot)=f(\cdot;\tau_t,\lambda_t)$,
$t\in \T$.

\vs3

\no For $\Ac=\M^d_+$, the conditions $\HFa$,  $\HFb$ and $\HL$
hold (this is a case where $-e_{ji} \notin \Ac$). However, by
construction, $\HNF$ does not hold except when
$N^0(F)=\{0\}$. As in  perfect information models,  this is the case if
$\lambda^{ij}+\lambda^{ji}>0$  for all $i,j\le d$.

\begin{Lemma} Fix $t\in \T$ and assume that for all $B\in \Hc_t$ there is $B'\subset B$
with positive probability such that, for all $i,j,k\le d$, $(1+\lambda^{ij}_t)$
$\le$ $(1+\lambda^{ik}_t)(1+\lambda^{kj}_t)$ and $\lambda^{ij}_t+\lambda^{ji}_t>0$ on $B'$.
Then    $N^0_t(F)=\{0\}$.
\end{Lemma}

\proof Fix $\eta\in L^0(\M^d_+;\Hc_t)$ and set $B:=\{\eta\ne 0\}$.  Under the above conditions,
on easily checks that the random cone
    \b*
    K_t&=&\{x\in \R^d~:~a\in \M^d_+,\; x+ \sum_{j\le d} a^{ji}
    -
    a^{ij}\tau^{ij}_t\left(1+\lambda^{ij}_t\right) \ge 0\;,\;\forall\;i\le d\}
    \e*
satisfies $K_t\cap(-K_t)$ $=$ $\{0\}$ on $B'$, see e.g. \cite{BT99} or \cite{KSR01}.
Since $N^0_t(F)$ $\subset$ $\{\xi \in L^0(\R^d;\Fc)~:~\xi \in
K_t\cap(-K_t)\}$, this shows that $F_t(\eta)=0$.
 \ep

\begin{Remark}{\rm The condition $\Pro{(1+\lambda^{ij}_t)
\le(1+\lambda^{ik}_t)(1+\lambda^{kj}_t)~|~\Hc_t}>0$ is natural
since otherwise it would be a.s. cheaper to transfer money from
$i$ to $j$ by passing through $k$ than directly. In this case, any
``optimal" strategy would induce an effective transaction cost
corresponding to $\tilde
\lambda^{ij}_t:=(1+\lambda^{ik}_t)(1+\lambda^{kj}_t)-1$. }
\end{Remark}

\no As argued in the introduction, if  $\tau$ or $\lambda$ are not
$\H$-adapted, transactions may   be non-reversible even when
transaction costs are equal to zero.

\begin{Lemma}  Assume that for some $t \in \T$ and $i\le d$,
$\tau^{ij}_t (1+\lambda^{ij}_t)$ is not $\Hc_t$-measurable for
all $j\le d$. Then, for all $\eta$ $\in L^0(\M^d_+;\Hc_t)$,
$F_t(\eta) \in N^0_t(F)$ implies $\sum_{j\le d}
\eta^{ji}+\eta^{ij}=0$.
\end{Lemma}

\proof  Fix $\eta \in L^0(\M^d_+;\Hc_t)$  such that $F_t(\eta)\in
N^0_t(F)$. By definition, there is $\tilde \eta$ $\in
L^0(\M^d_+;\Hc_t)$ such that $F_t(\eta)$ $=$ $-F_t(\tilde \eta)$.
Hence,
    \b*
    \sum_{j\le d}
    \eta^{ji} + \tilde \eta^{ji}
    &=&
    \sum_{j\le d}
    (\eta^{ij}+\tilde \eta^{ij})
    \tau^{ij}_t\;(1+\lambda^{ij}_t)\;.
    \e*
If $\sum_{j\le d} (\eta^{ij}+\tilde \eta^{ij})$ $\ne$ $0$, then the left hand-side term is $\Hc_t$-measurable while the right hand-side is not. It follows that both terms must be equal to $0$, so that $\sum_{j\le d} \eta^{ji}+\eta^{ij}=0$.
\ep\\

\no Since the conditions $\HFa$ and $\HFb$ hold, one can
argue  as in the above subsection to obtain  that,
if $N^0(F)=\{0\}$, then  $\NAr$ and   $\NAs$
are equivalent  to the existence of some $Z \in \Dc(F)$ which
must satisfy
    \b*
    \bar Z^i_t < \bar Z^j_t \;\bar \tau^{ji}_t \;(1+\bar \lambda^{ji}_t)
    \;\;\;,\;\;i,j\le d\;,\;t\in \T\;.
    \e*

\end{document}